\documentclass[11pt]{amsart}
\usepackage{amsmath,amssymb,amsthm, amscd}

\title{Symplectic torus bundles and group extensions}
\author{Peter J. Kahn}
\date{April 15, 2004}
\newcommand{\bundle}{\ensuremath{\xi: F
\overset{i}{\rightarrow}\  E\overset{p}{\rightarrow}\ B}}
\newcommand{\tbundle}{\ensuremath{\xi: T^{2}
\overset{i}{\rightarrow}\  E\overset{p}{\rightarrow}\ B}}
\newcommand{\bb}[1]{\ensuremath{\mathbb{#1}}}
\newcommand{\dd}{\ensuremath{d^{0,2}_{2}}}

\theoremstyle{definition}

\theoremstyle{remark}

\newtheorem*{remark}{Remark}
\newtheorem{example}{Example}

\theoremstyle{plain}
\newtheorem{theorem}{Theorem}[section]
\newtheorem*{ttheorem}{Theorem}

\newtheorem{Cor}[theorem]{Corollary}
\newtheorem{lemma}{Lemma}
\newtheorem{cor}[lemma]{Corollary}
\newtheorem{prop}{Proposition}[section]

\newtheorem{corr}[prop]{Corollary}

\begin{document}
\maketitle

\begin{abstract}
Symplectic torus bundles $\xi:T^{2}\rightarrow E\rightarrow B$ are classified by the second
cohomology group of $B$ with local coefficients $H_{1}(T^{2})$. For $B$ a compact, orientable
surface, the main theorem of this paper gives a necessary and sufficient condition on the
cohomology class corresponding to $\xi$ for $E$ to admit a symplectic structure compatible with
the symplectic bundle structure of $\xi$ : namely, that it be a torsion class. The proof is based
on a group-extension-theoretic construction of J. Huebschmann~\cite{jH}. A key ingredient is the
notion of fibrewise-localization.
\end{abstract}

\section{Introduction}\label{S:intro}

A symplectic $F$-bundle in this paper is a smooth fibre bundle $\bundle$ whose structure group is
the group of symplectomorphisms $Symp(F,\sigma)$ for some symplectic form $\sigma$ on $F$. For
such a bundle, the fibres $F_{b}=p^{-1}(b)$ admit canonical symplectic forms $\sigma_{b}$, the
pullbacks of $\sigma$ via symplectic trivializations.  A natural question to ask about $\xi$ is
under what conditions the forms $\sigma_{b}$ ``piece together'' to produce a symplectic form on
$E$. More exactly, when is there a \emph{closed}\  2-form $\beta$ on $E$ such that
\begin{equation}\label{E:main}
\beta|F_{b}=\sigma_{b},\quad\quad\ \mbox{for all}\  b\in\ B,
\end{equation}
with $\beta$ non-degenerate?  When $B$ is connected, an argument of W. Thurston (cf.~\cite[page
199]{MS}) shows that a closed 2-form $\beta$ satisfying (1) exists if and only if the de Rham
cohomology class of $\sigma$ is contained in $image(i^{\ast}:H^{2}_{DR}(E)\rightarrow\
H^{2}_{DR}(F))$. Thurston further shows that when such a $\beta$ exists and  $E$ is compact and
$B$ is symplectic, then $\beta$ may be modified to be non-degenerate while still satisfying (1).
McDuff and Salamon~\cite[page~202]{MS} use Thurston's result to settle the existence question for
a large family of surface bundles:

\begin{ttheorem}
Suppose that $F$ is a closed, oriented, connected surface of genus $\neq$ 1, and let
$\xi:F\rightarrow E\rightarrow B$ be a symplectic $F$-bundle with $B$ a compact, connected
symplectic manifold.  Then, $E$ admits a symplectic structure inducing the given structures on
the fibres. $\quad \square$
\end{ttheorem}
Their argument does not apply to the case of torus bundles, however; indeed, they present the
following simple counterexample in that case.  Consider the composition
\[   S^{1}\times S^{3} \overset{pr}{\rightarrow}\ S^{3} \overset{\mathfrak{H}}{\rightarrow}\ S^{2},\]
where $\mathfrak{H}$ is the well-known Hopf map. This composition is the projection of a
symplectic torus bundle.  No symplectic form can exist on the total space $S^{1}\times S^{3}$,
however, because $H^{2}_{DR}(S^{1}\times S^{3})=0$.

\subsection{The results}

This paper obtains a necessary and sufficient condition for the existence of $\beta$ in the case
of symplectic torus bundles over surfaces. Before stating our main result, however, we remind the
reader of some subsidiary facts.  For any fibre bundle $\xi:F\rightarrow\ E\rightarrow\ B$ with
group $G$, the action of $G$ on $F$ produces a $\pi_{0}(G)$-action on the homology and cohomology
of $F$.  When $B$ is a pointed space, there is a well-defined homomorphism
$\pi_{1}(B)\rightarrow\pi_{0}(G)$ that gives each homology or cohomology group of $F$ the
structure of a $\bb{Z}\pi_{1}(B)$-module.  Now suppose that $F$ is the $2$-torus $T^{2}$ and
$G=Symp(T^{2},\sigma)$.  It is not hard to show (see Appendices A and B) that $\pi_{0}(G)\approx
SL(2,\bb{Z})$ and that the $\pi_{0}(G)$-action on $H_{1}(T^{2})$ may be identified with the
natural action of $SL(2,\bb{Z})$ on $\bb{Z}^{2}$.  Given any representation
$\rho:\pi_{1}(B)\rightarrow\pi_{0}(G)=SL(2,\bb{Z})$, we let $\bb{Z}^{2}_{\rho}$ denote the
corresponding $\bb{Z}\pi_{1}(B)$-module.
\begin{prop}~\label{P:classification}
Assume that $B$ has the homotopy type of a pointed, path-connected $CW$ complex, and choose any
representation $\rho:\pi_{1}(B)\rightarrow\ SL(2,\bb{Z})$.  Then there is a natural, bijective
correspondence between the equivalence classes of symplectic torus bundles over $B$ inducing the
module structure $\bb{Z}^{2}_{\rho}$ on $H_{1}(T^{2})$ and the elements of
$H^{2}(B;\bb{Z}^{2}_{\rho})$, the second cohomology group of $B$ with local coefficients
$\bb{Z}^{2}_{\rho}$.\hfill$\square$
\end{prop}
\begin{remark}
We call the cohomology class corresponding to the symplectic torus bundle  $\xi$ the
\emph{characteristic class} of $\xi$ and denote it by $c(\xi)$. The characteristic class $c(\xi)$
vanishes if and only if $\xi$ admits a section.  When the representation $\rho$ is trivial,
$c(\xi)=0$ if and only if $\xi$ is trivial.
\end{remark}

The proposition and remark follow immediately from known, classical results of algebraic
topology, as  described in Appendices A and B.

\vspace{.3in}

We can now state the main result of this paper.

\begin{theorem}\label{T:main}
Suppose that $\xi$ is a symplectic torus bundle over a connected surface $B$.  Then the total
space of  $\xi$ admits a closed form $\beta$ satisfying (1) if and only if the characteristic
class $c(\xi)$ is a torsion element of $H^{2}(B;\bb{Z}^{2}_{\rho})$. If, in addition, $B$ is
compact and orientable and such a form exists, it can be chosen to be a symplectic form.
\end{theorem}

The last statement of the theorem is simply an application of Thurston's argument mentioned
above. So our proof of the theorem focuses exclusively on the existence of a closed $2$-form
$\beta$ satisfying~(\ref{E:main}).

\vspace{.2in}

The  following consequences of the theorem are almost immediate. We give proofs in
\S~\ref{S:Cor}.

\begin{Cor}\label{C:finite}
Let $B$ be a connected surface, and let $\rho:\pi_{1}(B)\rightarrow SL(2,\bb{Z})$ be a
representation.  Among the symplectic torus bundles over $B$ that induce the representation
$\rho$, there are, up to equivalence, only  finitely many  whose total spaces admit closed forms
$\beta$ satisfying~(\ref{E:main}).
\end{Cor}
\begin{Cor}\label{C:principal}
Every principal torus bundle has a canonical structure as a symplectic torus bundle. Let
$\xi:T^{2}\rightarrow E\rightarrow B$ be such a bundle, with $B$ a connected surface. Then, $E$
\emph{fails} to admit a closed $2$-form $\beta$ satisfying (1) if and only if $B$ is closed and
orientable \emph{and} $\xi$ is non-trivial.
\end{Cor}

A specialization of this corollary perhaps deserves a separate statement.

\begin{Cor}\label{C:torusaction}
Suppose the closed, connected symplectic $4$-manifold $E$ admits a free $T^{2}$-action such that
the orbits are symplectic submanifolds.  Then, as $T^{2}$-manifolds, $E\approx\ T^{2}\times\
(E/T^{2}). $
\end{Cor}
\begin{remark}
There does not appear to be a reasonable, \emph{non-trivial} sense in which the
$T^{2}$-equivariant diffeomorphism of this corollary can be taken to be a symplectomorphism.
There is simply too much leeway allowed by the hypotheses for symplectic forms on $E$.
\end{remark}

The proof of Theorem~\ref{T:main}\  breaks into three cases according as the base surface $B$ is
non-closed, closed of genus zero, and closed of genus different from zero. The first two cases
are substantially easier than the third and are proved at the end of this section and in \S\ 4,
respectively.  In these two cases, the theorem reduces to the following propositions.

\begin{prop}\label{P:mainnonclosed}
 Every symplectic torus bundle over a connected, non-closed surface admits a section and has a total space
that admits a closed $2$-form $\beta$ satisfying (1).
\end{prop}

\begin{prop}\label{P:maingenuszero}
    \begin{enumerate}
        \item The total space of a symplectic torus bundle over $S^{2}$ admits a closed $2$-form
        satisfying (1) if and only if the bundle is trivial.
        \item Let $E$ be the total space of a symplectic torus bundle $\xi$ over $\bb{R}P^{2}$.
        If the representation $\rho$ corresponding to $\xi$ is trivial,  then $E$ admits a
        closed $2$-form $\beta$ satisfying (1).  If $\rho$ is non-trivial, then $E$ admits
        such a $2$-form if and only if $c(\xi)=0$, that is, if and only if  $\xi$ admits a section.
    \end{enumerate}
\end{prop}

The case in which $B$ is a closed surface of genus $\neq 0$ forms the heart of the paper and
occupies \S\S2,3.  The following two examples suggest the variety of concrete possibilities in
this case. In both examples the base space $B$ is itself the torus $T^{2}$. Thus, in both, the
representation $\rho$ is a homomorphism $\pi_{1}(T^{2})=\bb{Z}^{2}\rightarrow\ SL(2, \bb{Z})$.

\begin{example}
For any $(a,b)\in \bb{Z}^{2}$,  define $\rho$ by the equation
\[
    \rho(a,b) =
    \begin{pmatrix}
      1 & b\\
      0 & 1
    \end{pmatrix}.
\]
In this example, one computes that the bundles are classified by
$H^{2}(T^{2};\mathbb{Z}^{2}_{\rho})=\mathbb{Z}.$ Consequently, up to equivalence, there is only
one torus bundle $\xi$---namely, the one satisfying $c(\xi)=0$--- for which the total space
admits a symplectic form satisfying (1). According to the classification, this is the unique
bundle admitting a section. The total space of $\xi$ is the renowned Kodaira-Thurston manifold,
the earliest known example of a symplectic manifold that is not K\"{a}hler (cf.
\cite[page~89]{MS}).
\end{example}

\begin{example}
Let $m$ and $n$ be any fixed integers $\geq 0$. Then, for $(a,b)\in \bb{Z}^{2}$, define $\rho$ by
\[
    \rho(a,b) =
    \begin{pmatrix}
      -2mn+1 & 2mn^{2}+n\\
      -m & mn+1
    \end{pmatrix}
^{a+b}.
\]
In this example, the bundles are classified by
$H^{2}(T^{2};\mathbb{Z}^{2}_{\rho})=\bb{Z}_{m}\oplus\mathbb{Z}_{n}$. So, when $m,n \neq 0$, there
are exactly $mn$ symplectic torus bundles over the torus, and, for every one of them, the total
space admits the desired symplectic form.
\end{example}

Both examples proceed by computing $H^{2}$ and then applying Theorem~\ref{T:main}. The
computation begins with Poincar\'{e} duality for $T^{2}$ (with twisted coefficients), which
implies that the desired result is just the group of coinvariants of the module
$\mathbb{Z}^{2}_{\rho}$ (cf. \cite[page~57]{kB}). We leave this computation to the reader.

\subsection{Reformulating Thurston's criterion}

We conclude this introduction with a brief reformulation of Thurston's cohomology criterion for
the existence of the desired closed $2$-forms $\beta$ in the context of symplectic torus bundles.
This will immediately imply Proposition~\ref{P:mainnonclosed}.

Thurston's criterion is stated in our opening paragraph in terms of de Rham cohomology, but
clearly, by de Rham's theorem, it may be equivalently stated in terms of singular cohomology with
real coefficients. In fact, a further easy reduction is desirable: namely, we pass to rational
coefficients. Indeed, note that since $H^{2}(T^{2};\bb{R})\approx \bb{R}$, the existence of a
non-trivial class in the image of $i^{\ast}:H^{2}(E;\bb{R})\rightarrow\ H^{2}(T^{2};\bb{R})$ is
equivalent to the surjectivity of this map, and this in turn  is easily checked to be equivalent
to the surjectivity of $i^{\ast}:H^{2}(E;\bb{Q})\rightarrow\ H^{2}(T^{2};\bb{Q})\approx\ \bb{Q}$.

Now using rational coefficients, we  consider the Serre cohomology spectral sequence for the
symplectic torus bundle \mbox{$\tbundle$}, for which the $E_{2}$-term is given by
\[E^{p,q}_{2}=H^{p}(B;H^{q}(T^{2};\bb{Q})).\]
Therefore, $E^{0,2}_{2}=H^{0}(B;H^{2}(T^{2};\bb{Q}))=H^{2}(T^{2};\bb{Q})^{\pi_{1}(B)}$,  the
group of $\pi_{1}(B)$-invariant classes in $H^{2}(T^{2};\bb{Q})$.  But $\pi_{1}(B)$ acts via
symplectomorphisms, which are orientation-preserving, so $E^{0,2}_{2}= H^{2}(T^{2};\bb{Q})$. Now
when $B$ is a surface, its cohomology vanishes above dimension two, so that $\dd$ is the only
possibly non-trivial differential issuing from $E^{0,2}_{r},  r\geq\ 2$. Thus
$ker(\dd:H^{2}(T^{2};\bb{Q})\rightarrow\ H^{2}(B; H^{1}(T^{2};\bb{Q})))= E^{0,2}_{\infty}$, which
equals $i^{\ast}(H^{2}(B;\bb{Q}))$.   Therefore, in this context, Thurston's cohomology criterion
becomes
\begin{equation}\label{E:differential}
\dd = 0.
\end{equation}

\noindent\emph{Proof of Proposition~\ref{P:mainnonclosed}}: Proposition~\ref{P:mainnonclosed} now
follows easily, using the fact that every connected, non-closed surface has the homotopy type of
a $1$-dimensional simplicial complex. Every $F$-bundle over such a base space admits a section
when $F$ is path-connected.  Moreover, the target of $\dd$, namely
$H^{2}(B;H^{1}(T^{2};\bb{Q}))$, is identically zero, so~(\ref{E:differential}) is
satisfied.\hfill$\square$

\vspace{.3in}

To conclude this introduction, I am pleased to to acknowledge my indebtedness to K. Brown for a
number of very helpful conversations during the preparation of this paper.

\section{An interpretation of the main theorem in terms of group extensions}\label{S:Extend}

Let $B$ be a connected, closed surface of genus $\neq\ 0$ and fundamental group $\pi$.  As is
well known, $B$ is a $K(\pi, 1)$, and so one sees easily that the homotopy exact sequence of the
symplectic torus bundle $\tbundle$ collapses to the short exact sequence
\begin{equation}\label{E:mainextension}
\bb{E}:\quad\quad \bb{Z}^{2}\overset{i_{\ast}}{\rightarrowtail}\
G\overset{p_{\ast}}{\twoheadrightarrow}\ \pi,
\end{equation}
which will be convenient to regard as a group extension of $\pi$ by $\bb{Z}^{2}$.  Thus, the
group $G$ equals $\pi_{1}(E)$, and $E$ is a $K(G,1)$.  Huebschmann~\cite{jH}\ uses the cohomology
spectral sequence of (3) (which is the same as the Serre spectral sequence of $\xi$) and obtains
group-extension-theoretic interpretations of some of its differentials.  We are interested in his
interpretation of
\[\dd: H^{2}(\bb{Z}^{2};\bb{Q})\rightarrow\ H^{2}(\pi;H^{1}(\bb{Z}^{2};\bb{Q}).\]
Here, we follow Huebschmann and use group-cohomology notation for the cohomology groups, but of
course these are the same as the cohomology groups of the base and fibre of $\xi$ as before.
Since $2$-dimensional group cohomology classifies group extensions with abelian kernel, the map
$\dd$ may be regarded as mapping extensions of $\bb{Z}^{2}$ by $\bb{Q}$---more precisely,
\emph{central} extensions, since $\bb{Z}^{2}$ acts trivially on $\bb{Q}$---to extensions of $\pi$
by $H^{1}(\bb{Z}^{2};\bb{Q})$. Huebschmann presents a construction that uses $\bb{E}$ to pass
from an extension $\bb{E}_{1}$ of the first kind to an extension $\bb{E}_{2}$ of the second.

\vspace{.2in}

\subsection{Huebschmann's construction}\label{s:Huebsch}

Let $\bb{E}_{1}$ denote an arbitrary \emph{central} extension of $\bb{Z}^{2}$ by $\bb{Q}$
\begin{equation}\label{E:E1}
\bb{Q}\rightarrowtail\ G_{1}\overset{r_{1}}{\twoheadrightarrow}\ \bb{Z}^{2}.
\end{equation}
We follow Huebschmann by using $\bb{E}$ and $\bb{E}_{1}$ to construct an extension $\bb{E}_{2}$
\begin{equation}~\label{E:extension2}
H^{1}(\bb{Z}^{2};\bb{Q})\rightarrowtail\ G_{2}\twoheadrightarrow\ \pi.
\end{equation}
We do this in several steps.

\vspace{.2in}

\noindent Step (a): Since $\bb{Z}^{2}$ is normal in $G$, inner automorphisms of $G$ determine
automorphisms of $\bb{Z}^{2}$, which give a representation
\begin{equation*}
\rho:\pi\rightarrow\ Aut(\bb{Z}^{2}) = GL(2,\bb{Z}).
\end{equation*}
We shall make use of the composition
\begin{equation*}
G\overset{p_{\ast}}{\twoheadrightarrow} \pi\overset{\rho}{\rightarrow} GL(2,\bb{Z}).
\end{equation*}
Further, every automorphism of $G_{1}$ determines a representation
\begin{equation}\label{E:rho1}
\rho_{1}:Aut(G_{1})\rightarrow\ Aut(G_{1}/\bb{Q}) = GL(2,\bb{Z}).
\end{equation}

The homomorphisms $\rho\circ  p_{\ast}$ and $\rho_{1}$ allow us to form the fibre product
$\Pi=G\times_{GL(2,\bb{Z})}Aut(G_{1})$. Let $p_{1}$ and $p_{2}$ denote the projections
$\Pi\rightarrow\ G$, $\Pi\rightarrow\ Aut(G_{1})$, respectively.

\vspace{.2in}

\noindent Step (b):  Combining (3) and (4), we have a composite homomorphism
\begin{equation}\label{E:lambda}
\lambda:G_{1}\overset{r_{1}}{\twoheadrightarrow}\ \bb{Z}^{2}\overset{i_{\ast}}{\rightarrowtail}\
G.
\end{equation}
and a homomorphism $\mu:G_{1}\rightarrow\ \Pi$ given by
\begin{equation}\label{E:mu}
\mu(x)=(\lambda(x),\iota_{x}),
\end{equation}
where $\iota_{x}$ denotes inner automorphism by $x$. It is not hard to check that
\begin{equation*}
\rho\circ p_{\ast}(\lambda(x))=\rho_{1}(\iota_{x})= I,
\end{equation*}
where $I$ is the $2\times\ 2$ identity matrix in $GL(2,\bb{Z})$.  Therefore, $\mu$ does indeed
take values in $\Pi$.  Let $G_{2}$ denote the quotient $\Pi/im(\mu)$ and $\lambda_{2}$ the
projection $\Pi\rightarrow\ G_{2}$.

\vspace{.2in}

\noindent Step(c): Note that $\mu$ vanishes on $ker(r_{1})$ so that it factors as
$G_{1}\overset{r_{1}}{\twoheadrightarrow} \bb{Z}^{2}\rightarrowtail\ \Pi$, where  the second map
lifts the injection $i_{\ast}:\bb{Z}^{2}\rightarrowtail\ G$.  It follows that $p_{1}$ maps
$im(\mu)$ bijectively onto $im(i_{\ast})$ , which implies that $p_{1}$ descends to a surjection
$r:G_{2}\twoheadrightarrow\ \pi$, and $\lambda_{2}$ maps $ker(p_{1})=H^{1}(\bb{Z}^{2};\bb{Q})$
isomorphically onto $ker(r)$. Therefore,  $r:G_{2}\twoheadrightarrow\pi$ is an extension of $\pi$
by $H^{1}(\bb{Z}^{2};\bb{Q})$, which is the desired extension $\bb{E}_{2}$ (see
(\ref{E:extension2}) above). The following diagram of exact sequences summarizes the situation:

\begin{equation*}
    \begin{CD}
    @.     0           @.            0                  @.             0\\
     @.     @VVV                    @VVV                            @VVV\\
    @.     \bb{Q}     @>0>>        H^{1}(\bb{Z}^{2};\bb{Q})    @=  H^{1}(\bb{Z}^{2};\bb{Q})\\
     @.      @VVV                    @VVV                            @VVV\\
    @.     G_{1}      @>\mu>>      \Pi       @>\lambda_{2}>>         G_{2} @>>> 0\\
     @.      @Vr_{1}VV                    @Vp_{1}VV                   @VrVV\\
    0 @>>> \bb{Z}^{2}  @>i_{\ast}>>        G       @>p_{\ast}>>                \pi @>>> 0\\
    @.     @VVV                    @VVV                         @VVV\\
    @.      0           @.            0          @.                     0
    \end{CD}
\end{equation*}

\vspace{.2in}

 Now let $c(\bb{E}_{1})$ and $c(\bb{E}_{2})$ denote the cohomology classes of the
extensions $\bb{E}_{1}$ and $\bb{E}_{2}$ , respectively.
\begin{ttheorem}[Huebschmann,~\cite{jH} ]
$\dd(c(\bb{E}_{1}))= c(\bb{E}_{2}).$
\end{ttheorem}

Huebschmann's result allows us to analyze properties of $\dd$ (e.g., condition (2)) by applying
his construction to a certain family of central extensions. Note, however, that the family we are
interested in may be described as $H^{2}(\bb{Z}^{2};\bb{Q})\approx \bb{Q}$, a $1$-dimensional
vector space over $\bb{Q}$.  So, to determine the vanishing of $\dd$, it suffices to analyze
Huebschmann's construction for any single central extension of $\bb{Z}^{2}$ by $\bb{Q}$ that
represents a non-zero element of cohomology.  We describe such an extension shortly, but first we
must make a short preparatory digression.

\subsection{Fibrewise-localization}

The theory of localization in algebraic topology has been well known since the work of Quillen,
Sullivan, Bousfield, Kan, Dwyer, Hilton, Mislin and others.  We summarize only that small
fragment of the subject that we need here.  A useful reference for the reader is~\cite{HM}.  We
shall confine ourselves to \emph{localizing at 0}, i.e., to rationalization, although most of
what we describe applies to the general case.

Localization of a nilpotent group $N$ is equivalent to localization of the Eilenberg-MacLane
space $K(N,1)$.  We'll use the language of groups here, however, rather than that of topology.
For the moment, we restrict entirely to nilpotent groups.  A \emph{local group} may be defined
here as a nilpotent group that is uniquely $p$-divisible for all primes $p$.  A localization of
the nilpotent group $N$ consists of a localization homomorphism (or localization map)
$\ell:N\rightarrow\ N_{0}$, where $N_{0}$ is local, such that $\ell$ is universal for
homomorphisms of $N$ into local groups (i.e., every such homomorphism $h:N\rightarrow\ L$ factors
as $h_{0}\ell$ for a unique homomorphism $h_{0}:N_{0}\rightarrow\ L$).  $N_{0}$ and $\ell$ are
uniquely determined up to the obvious equivalence.  When $N$ is abelian, $N_{0}$ may be taken to
be $N\otimes\bb{Q}$ and $\ell$ given by $x\mapsto\ x\otimes\ 1$.   A key fact about localization
is that localization maps induce localization homomorphisms of homology. Localization respects
exact sequences.  Indeed, it is not hard to show that, given any exact sequence $\bb{S}$ of
nilpotent groups, we may localize its terms and maps, obtaining an exact sequence $\bb{S}_{0}$ of
local groups and a map of exact sequences $\ell_{\bb{S}}:\bb{S}\rightarrow\bb{S}_{0}$ that
localizes the individual terms. Thus, we may apply this to group extensions in which all the
groups are nilpotent.

Let \[\bb{S}:   N^{\prime}\rightarrowtail\ N\twoheadrightarrow\ N^{\prime\prime}\] be a short
exact sequence of nilpotent groups, and let \[\bb{S}_{0}:   N^{\prime}_{0}\rightarrowtail\
N_{0}\twoheadrightarrow\ N^{\prime\prime}_{0}\] denote its localization. Then $\ell_{\bb{S}}$ may
be thought of as a triple of localization maps $(\ell_{N^{\prime}}, \ell_{N},
\ell_{N^{\prime\prime}})$. We use $\ell_{N^{\prime\prime}}:N^{\prime\prime}\rightarrow\
N^{\prime\prime}_{0}$ to pull back the sequence $\bb{S}_{0}$ to an exact sequence
\[\bb{S}_{f0}: N^{\prime}_{0}\rightarrowtail\ N_{f0}\twoheadrightarrow\ N^{\prime\prime},\] which we
call the \emph{fibrewise-localization} of $\bb{S}$. The pullback construction  produces a natural
map of exact sequences $\ell_{f}:\bb{S}\rightarrow \bb{S}_{f0}$ which on $N^{\prime\prime}$ is
just the identity and on $N^{\prime}$ is just the localization map
$\ell_{N^{\prime}}:N^{\prime}\rightarrow\ N^{\prime}_{0}$.

While this construction is perfectly valid, we want to use fibrewise-localization in the case of
group extensions with abelian kernel without assuming any nilpotency restrictions. So we present
another construction, valid for all such extensions.  Consider a group extension with abelian
kernel $A$,
\begin{equation}\label{E:extensionABC}
\bb{S}: A\rightarrowtail\ B\twoheadrightarrow\ C,
\end{equation}
and consider any normalized $2$-cocycle $\phi$ associated with $\bb{S}$.  This is a function
$\phi:C\times\ C\rightarrow\ A$ subject to normalization and $2$-cocyle identities
(cf.~\cite[pp.91 ff.]{kB}). $\phi$ is defined by choosing a \emph{function} $C\rightarrow\ B$
that splits the surjection $B\twoheadrightarrow\ C$ in (\ref{E:extensionABC}) and measuring how
far this deviates from being a homomorphism. Now, form the composite $C\times\
C\overset{\phi}{\rightarrow}\ A\overset{\ell}{\rightarrow}\ A_{0}$, where $\ell$ is a
localization map. This composite is a new normalized $2$-cocycle for an extension of $C$ by
$A_{0}$.  We define this extension to be the fibrewise-localization of $\bb{S}$ and denote it by
$\bb{S}_{f0}$.  There is an obvious map of extensions $\bb{S}\rightarrow\ \bb{S}_{f0}$ with the
same properties as before. It is not hard to show, using basic facts about extensions, that,  up
to equivalence of extensions, this construction is independent of the initial choice of
$2$-cocycle $\phi$ corresponding to $\bb{S}$ and independent of the choice of localization map
$\ell$, and it coincides with our earlier description of fibrewise-localization for nilpotent
extensions of nilpotent groups with abelian kernels.  Note also that this construction shows that
if $c(\bb{S})$ and $c(\bb{S}_{f0})$ are the cohomology classes of the corresponding extensions
(i.e., the cohomology classes of the corresponding $2$-cocycles), then the  homomorphism
$H^{2}(C;A)\rightarrow\ H^{2}(C;A_{0})$ induced by the localization map $\ell:A\rightarrow\
A_{0}$ sends $c(\bb{S})$ to $c(\bb{S}_{f0})$.

We now present a  useful and well-known extension of $\bb{Z}^{2}$ by $\bb{Z}$.

The discrete Heisenberg group $\mathcal{H}$ may be described as the set $\bb{Z}^{3}$ of all
integer triples with the following multiplication
\begin{equation}\label{E:Heisenbergop}
(a,b,c)\bullet(x,y,z)=(a+x+bz,b+y,c+z).
\end{equation}
The center $Z\lbrack\mathcal{H}\rbrack$ and commutator $\lbrack\mathcal{H},\mathcal{H}\rbrack$
both equal $\bb{Z}=\bb{Z}\times\ 0\times\ 0$, so that we clearly obtain the central extension
\begin{equation*}
\bb{H}:\quad\ \bb{Z}\rightarrowtail\  \mathcal{H}\twoheadrightarrow\ \bb{Z}^{2}.
\end{equation*}
We call this the \emph{Heisenberg extension}. The following result about $\bb{H}$ is well known.
For the convenience of the reader, we present a proof due to K. Brown.
\begin{lemma}
The cohomology class $c(\bb{H})$\  generates $H^{2}(\bb{Z}^{2};\bb{Z})\approx\bb{Z}$.
\end{lemma}
\begin{proof}
Let the group $H$ be given by the presentation $<x,y:[x,[x,y]], [y,[x,y]]>$.  If
$a,b\in\mathcal{H}$ are the triples $(0,1,0), (0,0,1)$, respectively, then it is not hard to
check that they generate $\mathcal{H}$, that $[a,b]= (1,0,0)$, and that, accordingly, $a$ and $b$
satisfy the relations for $x$ and $y$ in $H$ above.  Therefore, the rule $x\mapsto\ a, y\mapsto\
b$ well-defines a surjective homomorphism $f:H\rightarrow\mathcal{H}$.  We let the reader check
that this is injective as well.  Thus, $H\approx \mathcal{H}$, so that, given any group
$H^{\prime}$ and elements $c,d\in\ H^{\prime}$ satisfying the stated relations, there is a unique
homomorphism $\mathcal{H}\rightarrow\ H^{\prime}$ sending $a$ to $c$ and $b$ to $d$.

We apply this last fact  to an arbitrary central extension $\bb{M}: \bb{Z}\rightarrowtail\
M\twoheadrightarrow \bb{Z}^{2}$, choosing the elements $c,d\in\ M$ to be arbitrary lifts of
$(1,0), (0,1)\in\bb{Z}^{2}$, respectively.  Let $h:\mathcal{H}\rightarrow\ M$ be the
corresponding homomorphism. $h$ clearly induces a map of extensions $\bb{H}\rightarrow\bb{M}$
which is the identity on $\bb{Z}^{2}$ and is an endomorphism on $\bb{Z}$, say multiplication by
some integer $k$.  By tracing out the definition of the $2$-cocycle corresponding to an
extension, it is easy to check that $c(\bb{M})= kc(\bb{H})$.  Thus, $c(\bb{H})$ generates
$H^{2}(\bb{Z}^{2};\bb{Z})\approx \bb{Z}$.
\end{proof}

We now define the extension of $\bb{Z}^{2}$ by $\bb{Q}$ that interests us: namely, it is the
fibrewise-localization of the Heisenberg extension, $\bb{H}_{f0}$.

\begin{cor}\label{C:basisclass}
$c(\bb{H}_{f0})$ is a basis element of the $1$-dimensional $\bb{Q}$ vector space
$H^{2}(\bb{Z}^{2};\bb{Q})$.
\end{cor}
\begin{proof}
Let $\ell_{\ast}:H^{2}(\bb{Z}^{2};\bb{Z})\rightarrow\ H^{2}(\bb{Z}^{2};\bb{Q})$ denote the
homomorphism induced by the coefficient injection $\bb{Z}\rightarrow\bb{Q}$.  As already
observed, $\ell_{\ast}$ maps $c(\bb{H})$ to $c(\bb{H}_{f0})$.  At the same time, it is clear that
$\ell_{\ast}$ is a localization map, essentially the same as the standard injection
$\bb{Z}\rightarrow\ \bb{Q}$.  Therefore, by the foregoing lemma, $c(\bb{H}_{f0})\neq\  0$, as
desired.
\end{proof}

\subsection{Reinterpreting the main theorem}

Let us return to the context with which this section opened: namely, to the symplectic torus
bundle $\tbundle$ with $B$  a closed, connected $K(\pi, 1)$ surface. The group $\pi$ acts via
symplectomorphisms on $H_{1}(T^{2}) = \bb{Z}^{2}$.  Thus, we have a representation $\rho$ and
corresponding (left) $\bb{Z}\lbrack\pi\rbrack$-module $\bb{Z}^{2}_{\rho}$, as explained before.
In a similar way, the cohomology group $H^{1}(T^{2};\bb{Q})\approx \bb{Q}^{2}$ receives the
structure of a $\bb{Z}\lbrack\pi\rbrack$-module.  We want this to be a left
$\bb{Z}\lbrack\pi\rbrack$-module also despite the contravariance of cohomology, so we use the
standard convention for this, which we may describe here as follows: Identify
$H^{1}(T^{2};\bb{Q})$ with $Hom(H_{1}(T^{2}),\bb{Q})$, and for any $\alpha\in\pi$, $h\in
Hom(H_{1}(T^{2}),\bb{Q})$, and $x\in\ H_{1}(T^{2})$, let $(\alpha h)(x)=h(\alpha^{-1}x)$.

 We now return to our use of group cohomology notation in the following lemma, the proof of which
 is given in the next section.

\begin{lemma}
Let $D:H^{1}(\bb{Z}^{2};\bb{Q})\rightarrow\ H_{1}(\bb{Z}^{2};\bb{Q})$ denote Poincar\'{e}
duality, and let $\psi$ be the composite
\[\bb{Z}^{2}=H_{1}(\bb{Z}^{2};\bb{Z})\overset{\ell}{\rightarrow}\
H_{1}(\bb{Z}^{2};\bb{Q})\overset{D^{-1}}{\rightarrow}\ H^{1}(\bb{Z}^{2};\bb{Q}),\] where, here,
$\ell$  is the localization map induced by the usual injection $\bb{Z}\rightarrowtail\bb{Q}$.
Then, using the module structures described above, $\psi$ is a
$\bb{Z}\lbrack\pi\rbrack$-injection and a localization map.

Therefore,

\begin{equation*}
\psi_{\sharp}:H^{2}(\pi;\bb{Z}^{2})\rightarrow\ H^{2}(\pi;H^{1}(\bb{Z}^{2};\bb{Q}))
\end{equation*}
induced by $\psi$ is also a localization map.
\end{lemma}

We can now state a reinterpretation of Theorem~\ref{T:main} in this group-extension context.

\begin{theorem}\label{T:mainequation}
Let $\bb{H}_{f0}$ be the fibrewise-localization of the Heisenberg extension, and let $\bb{E}$ be
the group extension (~\ref{E:mainextension}) described at the start of \S 2.  Apply Huebschmann's
construction to these, obtaining an extension $\bb{E}_{2}$ as in (~\ref{E:extension2}). Then,
\[\psi_{\sharp}(c(\bb{E}))= - c(\bb{E}_{2}).\]
\end{theorem}

 We prove Theorem~\ref{T:mainequation} in the next section. We close this section by using it
to prove Theorem~\ref{T:main} in case $B$ is closed, connected of genus $\neq\ 0$:

\begin{proof}
Let $\tbundle$ be a symplectic torus bundle with corresponding group extension $\bb{E}$. As
discussed in Appendix C, the classes $c(\xi)$and $c(\bb{E})$ are the same, so we may deal
exclusively with the latter. Suppose it has finite order. Then, by Huebschmann's theorem and
Theorem~\ref{T:mainequation},
 \[\dd(c(\bb{H}_{f0}))= c(\bb{E}_{2})= -\psi_{\sharp}(c(\bb{E}))= 0.\]  By Corollary~\ref{C:basisclass} of \S 2.2,
 this implies that $\dd=0$, which is condition (\ref{E:differential}). Therefore, as already argued, the desired
 form $\beta$ exists. The converse follows by reversing the steps.
\end{proof}

\section{Proof of Theorem 2.1}

The basic idea of the proof of Theorem 2.1 is to produce suitable $2$-cocycles $f$ and $F$ for
the extensions $\bb{E}$ and $\bb{E}_{2}$, respectively, and then to show that, if $\psi_{\flat}$
is the chain map induced by $\psi$, then $\psi_{\flat}(f)=-F$. To carry this out, we need to be
more explicit about $\psi$ and about the groups and maps occurring in Huebschmann's construction.

\subsection{The map $\psi$}

We begin with a proof of Lemma 3 of \S2.

\begin{proof}
That $\psi=D^{-1}\ell$ is a localization map and injective is obvious.  Choose any $\alpha\in
\pi$, and let $a$ be a symplectomorphism of $T^{2}$ representing $\alpha$.  This is a degree-one
map. Therefore, the standard cap product identity yields $a_{\ast}Da^{\ast}=D$, or
$a_{\ast}D=D(a^{\ast})^{-1}$, that is $\alpha\ D=D\alpha$.  So, D is
$\bb{Z}\lbrack\pi\rbrack$-equivariant. That $\ell$ is also equivariant is immediate from
definitions.  Hence $\psi$ is a map of $\bb{Z}\lbrack\pi\rbrack$-modules.

It remains to show that $\psi_{\sharp}:H^{2}(\pi;\bb{Z}^{2})\rightarrow\
H^{2}(\pi;H^{1}(\bb{Z}^{2};\bb{Q}))$ is a localization map.  By definition, $\psi_{\sharp}$
factors as
\begin{equation*}
H^{2}(\pi;\bb{Z}^{2})\overset{\ell_{\sharp}}{\rightarrow} H^{2}(\pi;
\bb{Q}^{2})\overset{(D^{-1})_{\sharp}}{\underset{\approx}{\longrightarrow}} H^{2}(\pi;
H^{1}(\bb{Z}^{2};\bb{Q})).
\end{equation*}
So, $\psi_{\sharp}$ is equivalent to $\ell_{\sharp}$.  But $\pi$ is finitely-presented, hence of
type $FP_{2}$ ~(\cite[page 197]{kB}). It follows without difficulty that $\ell_{\sharp}$ is
equivalent to the standard localization map \mbox{$H^{2}(\pi;\bb{Z}^{2})\rightarrow
H^{2}(\pi;\bb{Z}^{2})\otimes \bb{Q}$}.
\end{proof}

For computations which follow below, it will be useful to obtain an alternative description of
$\psi$.  Accordingly, we let $e_{1}$ and $e_{2}$ be the standard generators of
$H_{1}(\bb{Z}^{2};\bb{Z})=\bb{Z}^{2}$; we may write $a_{1}e_{1}+a_{2}e_{2}$ as $(a_{1}, a_{2})$.
Let $e_{1}^{\ast}, e_{2}^{\ast}$ denote the basis of $H^{1}(\bb{Z}^{2};\bb{Q})$ dual to
$\ell(e_{1}), \ell(e_{2})$, using this to write elements of $H^{1}(\bb{Z}^{2};\bb{Q})$ as pairs.
Then, one easily computes, $\psi(e_{1})= e_{2}^{\ast}$ and $\psi(e_{2})=-e_{1}^{\ast}$, so that,
in pair notation,
\begin{equation}\label{E:psi}
\psi(a_{1},a_{2})= (-a_{2},a_{1}).
\end{equation}

\subsection{$\bb{E}$ and the $2$-cocycle $f$}

Recall that $\bb{E}$ is the extension
\begin{equation*}
\bb{Z}^{2}\ \overset{i_{\ast}}{\rightarrowtail}\ G\ \overset{p_{\ast}}{\twoheadrightarrow}\ \pi.
\end{equation*}
Choose an
 an arbitrary \emph{function} $s:\pi\rightarrow G$ splitting $p_{\ast}$ and define the normalized
 $2$-cocycle $f$ by the usual rule
\begin{equation}\label{E:cocycle}
i_{\ast}(f(x,y))= s(x)s(y)s(xy)^{-1}.
\end{equation}
Now $f$, together with the representation $\rho:\bb{Z}^{2}\rightarrow GL(2,\bb{Z})$ induced by
$\bb{E}$, can be used to form another extension $\bb{E}^{\prime}$ of $\pi$ as follows:  In the
cartesian product $\bb{Z}^{2}\times\pi$ define a group multiplication $\bullet$ by the rule
\begin{equation}\label{E:opinG}
(u,x)\bullet(v,y)= (u+\rho(x)(v)+f(x,y), xy).
\end{equation}
Define homomorphisms $\bb{Z}^{2}\rightarrowtail \bb{Z}^{2}\times\pi$ and
$\bb{Z}^{2}\times\pi\twoheadrightarrow\pi$  by the rules $u\mapsto(u,\epsilon)$ and
$(u,x)\mapsto\ x$, respectively, where $\epsilon$ denotes the identity of $\pi$.  These piece
together to give the extension $\bb{E}^{\prime}$.  It is a classical fact that $\bb{E}$ and
$\bb{E}^{\prime}$ are equivalent extensions, and so $c(\bb{E})=c(\bb{E}^{\prime})$.  Therefore,
without losing generality, we may assume that $\bb{E}=\bb{E}^{\prime}$.

 With this assumption, the map $\lambda:\mathcal{H}_{f0}  = G_{1}\rightarrow\ G$ defined in (\ref{E:lambda}) can now be
 expressed as follows:
\[\lambda(a,b,c)= (b,c,\epsilon),\]
where we omit extra parentheses when harmless.  We want to get a similar explicit representation
of the map $\mu$ used above to define $G_{2}$, and for this, we need some computational
information about $\mathcal{H}_{f0}$ and $Aut(\mathcal{H}_{f0})$.

\subsection{Computational information about $\mathcal{H}_{f0}$ and $Aut(\mathcal{H}_{f0})$}

    We shall always regard $\mathcal{H}$ as embedded
in $\mathcal{H}_{f0}$ via the inclusion $\bb{Z}^{3}\subseteq \bb{Q}\times\bb{Z}^{2}$.

Given elements $x$ and $y$ in some group, we let $^{x}y$ denote the conjugate $xyx^{-1}$.  The
following lemma may be easily derived by the reader from the definition of the operation
(\ref{E:Heisenbergop}).

\begin{lemma}
In $\mathcal{H}_{f0},\quad ^{(a,b,c)}(x,y,z)= (x+bz-cy, y, z),$ and
$[(a,b,c),(x,y,z)]=(bz-yc,0,0)$. $\square$
\end{lemma}
\begin{cor}
The center $Z[\mathcal{H}_{f0}]$ equals $\bb{Q}\times0\times0$, setwise and as abelian
groups.$\quad \square$
\end{cor}

Thus, the surjection $\mathcal{H}_{f0}\rightarrow\bb{Z}^{2}$ in $\bb{H}_{f0}$ is just the
projection $\mathcal{H}_{f0}\rightarrow \mathcal{H}_{f0}/Z[\mathcal{H}_{f0}]$.  Recall that we
have denoted this $r_{1}$ in our description of Huebschmann's construction (cf. (\ref{E:E1})).

\begin{lemma}
Every endomorphism $h$ of $\mathcal{H}$ (resp., $\mathcal{H}_{f0}$) is uniquely determined by the
values $h(0,1,0)$ and $h(0,0,1)$.
\end{lemma}
\begin{proof}
The result is obvious for $\mathcal{H}$, since $(0,1,0)$ and $(0,0,1)$ generate it.  So, suppose
$h$ is an endomorphism of $\mathcal{H}_{f0}$.  For the reason just given, $h|\mathcal{H}$ is
uniquely determined by the given values.  Assume for the moment $h$ takes $Z[\mathcal{H}_{f0}]$
to itself.  That is, the restriction of $h$ to the center may be identified with an endomorphism
of $\bb{Q}$. But every such endomorphism is uniquely determined by its value at any single
non-zero element. Therefore, $h|Z[\mathcal{H}_{f0}]$ is determined by
$h(1,0,0)=[h(0,1,0),h(0,0,1)]$. Since $\mathcal{H}_{f0}$ is generated by $\mathcal{H}\cup
Z[\mathcal{H}_{f0}]$, the result holds for $\mathcal{H}_{f0}$.

It remains to show that $h$ maps $Z[\mathcal{H}_{f0}]$ to itself.  By Corollary 5, every element
$z$ in the center is $q$-divisible for every prime $q$.  Therefore, the same holds for any
homomorphic image of $z$, for example, for $r_{1}(h(z))\in \bb{Z}^{2}$.  But the only element of
$\bb{Z}^{2}$ with this divisibility property is $0$.  So, $h(z)\in
ker(r_{1})=Z[\mathcal{H}_{f0}]$, as required.
\end{proof}

\begin{lemma}
For any triples $(a,b,c), (d,e,f)\in\ \mathcal{H}_{f0}$, there exists an endomorphism $h$ of
$\mathcal{H}_{f0}$ satisfying $h(0,1,0)=(a,b,c)$ and $h(0,0,1)=(d,e,f)$.  $h$ is an automorphism
if and only if the determinant
\begin{equation*}
    \begin{vmatrix}
      b & c\\
      e & f
    \end{vmatrix}= \pm 1
\end{equation*}
\end{lemma}
\begin{proof}
By Lemma 4 and Corollary 5, the commutator $[(a,b,c),(d,e,f)]$ belongs to $Z[\mathcal{H}_{f0}]$,
so by the argument in the proof of Lemma 1 of \S2.2, there is a unique homomorphism
$k:\mathcal{H}\rightarrow\ \mathcal{H}_{f0}$ satisfying $k(0,1,0)= (a,b,c)$ and
$k(0,0,1)=(d,e,f)$.  By Lemma 4, $k(1,0,0)=(bf-ec,0,0)$, so it belongs to $Z[\mathcal{H}_{f0}]$,
and there is a unique extension of $k|Z[\mathcal{H}]$ to an endomorphism of
$Z[\mathcal{H}_{f0}]$. Every element $y$ of $\mathcal{H}_{f0}$ can be written as a product $zx$,
with $z\in\ Z[\mathcal{H}_{f0}]$ and $x\in\ \mathcal{H}$, so we attempt to define $h$ by the
rule, $h(y)=k(z)k(x)$.  It is an easy exercise to verify that this gives a well-defined
endomorphism.  Now suppose that $h$ is an automorphism.  Then it induces an automorphism of
$\bb{Z}^{2}$ given by the matrix
\begin{equation*}
    \begin{pmatrix}
      b & c\\
      e & f
    \end{pmatrix},
\end{equation*}
which immediately shows that the stated determinant must equal $\pm 1$.  Conversely, if the
determinant is $\pm 1$, then by what was just said, the endomorphism of $\bb{Z}^{2}$ induced by
$h$ is an automorphism, and, by the equation $h(1,0,0)=(bf-ec,0,0)$, so is the endomorphism of
$Z[\mathcal{H}_{f0}]$.  The Five-Lemma then implies that $h$ is an automorphism.
\end{proof}

We now introduce some convenient `matrix' notation for automorphisms $h\in\
Aut(\mathcal{H}_{f0}).$  If $h(0,1,0)=(a,b,c)$ and $h(0,0,1)=(d,e,f)$, as above, we associate
with $h$ the matrix
\begin{equation*}
    \begin{pmatrix}
        a & d\\
        b & e\\
        c & f
    \end{pmatrix}.
\end{equation*}
We may occasionally wish to abbreviate this by letting, say, $u$ denote the top row and, say, $M$
the remaining $2\times 2$ submatrix and writing the above matrix  as
\begin{equation*}
    \begin{pmatrix}
        u\\
        M
    \end{pmatrix}.
\end{equation*}
Of course, the identity automorphism has the obvious matrix representation
\begin{equation*}
    \begin{pmatrix}
        0 & 0\\
        1 & 0\\
        0 & 1
    \end{pmatrix}.
\end{equation*}
Slightly less obvious, but useful, is the matrix representation of the inner automorphism
$\iota_{x}$, where $x=(a,b,c)$. An easy application of Lemma 4 and equation (\ref{E:psi}) above
shows that this is:
\begin{equation*}
    \begin{pmatrix}
        -c & b\\
        1 &  0\\
        0 &  1
    \end{pmatrix} =
    \begin{pmatrix}
        \psi(b,c)\\
        I
    \end{pmatrix},
\end{equation*}
where $I$ is the $2\times 2$ identity matrix. It is possible to work out the multiplication,
i.e., composition, in $Aut(\mathcal{H}_{f0})$ in terms of this notation, but the formula is
complicated and not particularly useful here---in addition to the usual quadratic terms of linear
algebra, there are also third and fourth order terms. We do record one special case, however:
namely, the case of elements of the kernel of the natural projection
$\rho_{1}:Aut(\mathcal{H}_{f0})\rightarrow\ GL(2,\bb{Z})$ in (\ref{E:rho1}).  In matrix notation,
these elements consist of all matrices of the form,
\begin{equation*}
    \begin{pmatrix}
        u\\
        I
    \end{pmatrix}.
\end{equation*}
  In this case, one computes easily that
\begin{equation*}
    \begin{pmatrix}
        u\\
        I
    \end{pmatrix}\circ
    \begin{pmatrix}
        v\\
        I
    \end{pmatrix}=
    \begin{pmatrix}
        u+v\\
        I
    \end{pmatrix}.
\end{equation*}
Thus, the kernel is isomorphic, as an abelian group, to $\bb{Q}^{2}$.  Now, in fact, we know this
for other reasons: the kernel is known to be isomorphic to  $Hom(\bb{Z}^{2},\bb{Q})\approx
H^{1}(\bb{Z}^{2};\bb{Q})\approx\bb{Q}^{2}$. However, it is convenient for our computations to
have an explicit realization as $\bb{Q}^{2}$.

The following lemma provides a critical ingredient in the proof of Theorem~\ref{T:mainequation}
and explains our use of fibrewise-localization:

\begin{lemma}
$\rho_{1}:Aut(\mathcal{H}_{f0})\rightarrow\ Aut(\bb{Z}^{2})=GL(2,\bb{Z})$ is a split surjection.
\end{lemma}
\begin{proof}
That $\rho_{1}$ is surjective is an immediate corollary of Lemma 7.  To show that it splits, we
consider the extension $H^{1}(\bb{Z}^{2};\bb{Q})\rightarrowtail\ Aut(\mathcal{H}_{f0})
\overset{\rho_{1}}{\twoheadrightarrow}\ GL(2,\bb{Z})$, which represents an element of
$H^{2}(GL(2,\bb{Z});H^{1}(\bb{Z}^{2};\bb{Q}))$.  Now, the virtual cohomological dimension of
$SL(2,\bb{Z})$ is $1$,~\cite[page~229]{kB}. That is, it possesses a finite-index subgroup of
cohomological dimension $1$.  Therefore, the same holds for $GL(2,\bb{Z})$.  It follows easily
that $H^{i}(GL(2,\bb{Z});V)=0$ for all $i\geq 2$ and all $\bb{Q}[GL(2,\bb{Z})]$-modules $V$.
Thus, $H^{2}(GL(2,\bb{Z});H^{1}(\bb{Z}^{2};\bb{Q}))=0$, implying that $\rho_{1}$ splits.
\end{proof}

\vspace{.2in}

Choose and fix an arbitrary (homomorphic!) splitting $\tau:GL(2,\bb{Z})\rightarrow\
Aut(\mathcal{H}_{f0})$.

\subsection{The proof of Theorem 2.1}

We begin by rewriting the definition of the map $\mu:\mathcal{H}_{f0}\rightarrow\
\Pi=G\times_{GL(2,\bb{Z})}Aut(\mathcal{H}_{f0})$ in terms of the notation just introduced. Recall
that, for $z\in \mathcal{H}_{f0}$, $\mu(z)=(\lambda(z), \iota_{z})$, as above in (\ref{E:mu}) and
ff. Setting $z=(a,b,c)$ and using results in \S\S\ 3.2, 3.3, we have
\begin{equation}\label{E:newmu}
\mu(a,b,c)=((b,c,\epsilon), \begin{pmatrix}
                -c & b\\
                 1 & 0\\
                 0 & 1
                  \end{pmatrix}).
\end{equation}

We now proceed to define a $2$-cocycle $F$ for the extension $\bb{E}_{2}$ by first defining a
function $t:\pi\rightarrow G_{2}$ that splits the surjection $r:G_{2}\twoheadrightarrow \pi$ .
Recall that the standard projection $\Pi\rightarrow  G_{2}=\Pi/im(\mu)$ is denoted $\lambda_{2}$.
For any $w\in\ \Pi$, let us write $\lambda_{2}(w)=[w]$.  Then, for any $x\in \pi$, we define
$t(x)$ by
\begin{equation}\label{E:t}
t(x)=[(0,0,x), \tau(\rho(x))].
\end{equation}
Now we define F by the usual formula:
\begin{equation}\label{E:jayef}
j(F(x,y))=t(x)t(y)t(xy)^{-1},
\end{equation}
where $j:H^{1}(\bb{Z}^{2};\bb{Q})\rightarrow G_{2}$ is the inclusion onto  $ker(r)$. Let us make
$j$ more explicit.  Choose any $\phi\in H^{1}(\bb{Z}^{2};\bb{Q})=Hom(\bb{Z}^{2},\bb{Q})$. Then
$j(\phi)$ is precisely the image under $\lambda_{2}$ of the following pair in
$G\times_{GL(2,\bb{Z})}Aut(\mathcal{H}_{f0})=\Pi$:
\begin{equation}\label{E:jay}
    ((0,0,\epsilon),\begin{pmatrix}
                \phi(e_{1}) & \phi(e_{2})\\
                1      &   0\\
                0      &   1
            \end{pmatrix}),
\end{equation}
where, as before, $e_{1}, e_{2}$ are the standard generators of $\bb{Z}^{2}$. Now using equation
(\ref{E:opinG}), which gives the multiplication in $G$, we can compute $t(x)t(y)$:
\begin{align}
t(x)t(y)&= [(0,0,x)(0,0,y), \tau(\rho(x))\tau(\rho(y))]\notag\\
    &=[(f(x,y),xy), \tau(\rho(xy))]\notag\\
    &=[(f(x,y),\epsilon), \begin{pmatrix}
                0\\
                I
                 \end{pmatrix}][(0,0,xy), \tau(\rho(xy))].\notag
\end{align}
Note that the second and third equalities follow from the definition of the multiplication in
$G$, as given in equation (\ref{E:opinG}), as well as the fact that $\tau$ and $\rho$ are
homomorphisms! Now, using equation (\ref{E:t}), we get
\begin{equation*}
    t(x)t(y)=[(f(x,y),\epsilon), \begin{pmatrix}
                0\\
                I
                  \end{pmatrix}]t(xy),
\end{equation*}
which, when combined with ({\ref{E:jayef}), yields

\begin{equation*}
j(F(x,y))= [(f(x,y),\epsilon), \begin{pmatrix}
                0\\
                I
                  \end{pmatrix}].
\end{equation*}
Setting $f(x,y)=(f_{1},f_{2})=f_{1}e_{1}+f_{2}e_{2}\in \bb{Z}^{2}$ and applying equations
(\ref{E:newmu}) and (\ref{E:jay}), this becomes
\begin{align}
j(F(x,y))&=[(0,0,\epsilon), \begin{pmatrix}
                f_{2} & -f_{1}\\
                1        &     0\\
                0        &     1
               \end{pmatrix}]\notag \\
    &=j(-\psi(f(x,y)))\notag.
\end{align}

Since $j$ is injective,  $\psi(f(x,y))= -F(x,y)$, or $\psi_{\flat}(f)=-F$. This completes our
proof of  Theorem 2.1\hspace{1in} $\square$

\vspace{.3in}

\section{The Main Theorem when $B=S^{2}$ or $\bb{R}P^{2}$}

Let $\xi: T^{2} \overset{i}{\rightarrow} E \overset{p}{\rightarrow}B$ be a symplectic torus
bundle with $B$ a closed genus zero surface.  In this case, Theorem~\ref{T:main} reduces to
Proposition~\ref{P:maingenuszero}, which we prove in this section by methods essentially
unrelated to our earlier arguments.

First we deal with the case $B=S^{2}$.

\noindent\emph{Proof of Proposition~\ref{P:maingenuszero}(a)} \nopagebreak[4]

As we explain in  Appendix B, the classification of symplectic torus bundles over a
simply-connected space is the same as the classification of principal torus bundles over that
space. It is well-known that when $B=S^{2}$ these are classified by $\pi_{1}(T^{2})$.  Indeed,
the homotopy class corresponding to $\xi$ may be described as follows (cf.~\cite[page 98]{nS}).
Consider the following portion of the exact homotopy sequence of $\xi$:
\[  \pi_{2}(S^{2}) \overset{\partial}{\rightarrow}\pi_{1}(T^{2})\rightarrow\pi_{1}(E)\rightarrow\ 0.\]
Then the required homotopy class is $\pm\partial(\iota)\in\pi_{1}(T^{2})$, where $\iota$ is the
class of the identity map of $S^{2}$.  Since $\pi_{1}(E)$ is a homomorphic image of
$\pi_{1}(T^{2})$, it is abelian and thus equals $H_{1}(E)$.  It follows that this last has rank
one or two according as $\xi$ is non-trivial or trivial, respectively.  By Poincar\'{e} duality,
which applies because $E$ is closed and orientable, the same is true of the rank of  $H^{3}(E)$.

We now turn to the following portion of the Wang sequence for $\xi$:
\[H^{2}(E;\bb{Q})\overset{i^{\ast}}{\rightarrow}H^{2}(T^{2};\bb{Q})\overset{\theta}{\rightarrow}H^{1}(T^{2};\bb{Q})\rightarrow H^{3}(E;\bb{Q})\rightarrow\ 0.\]
Clearly, $i^{\ast}$ in this sequence is onto when $H^{3}(E)$ has rank two and $0$ when $H^{3}(E)$
has rank one. Since the surjectivity of $i^{\ast}$ with rational coefficients is equivalent to
the existence of the desired form $\beta$, this concludes the proof of
Proposition~\ref{P:maingenuszero}(a). $\quad \square$

We now deal with the case $B=\bb{R}P^{2}$.  Let $\pi:S^{2}\rightarrow \bb{R}P^{2}$ be the double
cover, and let $\tilde{E}$ be the total space of the pullback $\pi^{\ast}\xi$, a symplectic torus
bundle over $S^{2}$.  Then we have the following lemma.

\begin{lemma}~\label{L:pullback}
The total space $E$ of $\xi$ admits a closed $2$-form $\beta$ satisfying (\ref{E:main}) if and
only if $\tilde{E}$ does.
\end{lemma}
\begin{proof}
\noindent(a)$\Rightarrow$\  Let $\bar{\pi}:\tilde{E}\rightarrow\ E$ be the bundle map over $\pi$
given by the pullback construction.  If $\beta$ is a closed $2$-form on $E$ satisfying (1), then
$\bar{\pi}^{\ast}(\beta)$ is a closed $2$-form on $\tilde{E}$ satisfying (1).

\noindent(b)$\Leftarrow$\  Let $b:\tilde{E}\rightarrow\tilde{E}$ be the non-trivial deck
transformation. It is not hard to check, using the definition of the pullback construction, that
$b$ maps fibres of $\tilde{E}$ to fibres so as to preserve the pullback symplectic structures.
Now let $\gamma$ be a closed $2$-form on $\tilde{E}$ satisfying (1), and define
\[\tilde{\beta}=\frac{1}{2}(\gamma+b^{\ast}\gamma).\] Since $\tilde{\beta}$ is invariant under
deck transformations it descends to a closed $2$-form $\beta$ on $E$.  It clearly also satisfies
(1), which implies the same for  $\beta$.
\end{proof}

This lemma immediately implies the first statement of Proposition~\ref{P:maingenuszero}(b).

\begin{cor}
Suppose that the representation $\rho:\pi_{1}(\bb{R}P^{2})\rightarrow GL(2,\bb{Z})$ is trivial.
Then $E$ admits a closed $2$-form $\beta$ satisfying (1).
\end{cor}
\begin{proof}
If the module structure on $\bb{Z}^{2}$ is trivial, then
$H^{2}(\bb{R}P^{2};\bb{Z}^{2}_{\rho})\approx (\bb{Z}_{2})^{2}$.  Clearly then the map
$\pi^{\ast}: H^{2}(\bb{R}P^{2};\bb{Z}^{2}_{\rho})\rightarrow H^{2}(S^{2};\bb{Z}^{2})\approx
\bb{Z}^{2}$ is trivial.  By the classification theorem, it follows that the pullback
$\pi^{\ast}(\xi)$ is trivial.  But Proposition 3(a) then implies that the total space of this
pullback admits the desired $2$-form.  Therefore, by the lemma, so does $E$.
\end{proof}

It remains to deal with the case $B=\bb{R}P^{2}$, $\rho$ non-trivial. Since we are dealing with a
symplectic torus bundle, $\rho$ must take values in $SL(2,\bb{Z})$, which easily implies that
$im(\rho)=\{\pm I\}$. We now consider the cohomology Serre spectral sequence of the covering
$\pi:S^{2}\rightarrow\bb{R}P^{2}$, which has
\[E^{p,q}_{2}=H^{p}(\bb{Z}_{2};H^{q}(S^{2};\bb{Z}^{2}_{\rho}))\]
and converges to $H^{\ast}(\bb{R}P^{2};\bb{Z}^{2}_{\rho})$.  Here, the group
$H^{q}(S^{2};\bb{Z}^{2}_{\rho})$ is the ordinary cohomology of $S^{2}$ with $\bb{Z}^{2}$
coefficients, but the action of $\bb{Z}_{2}$ is a joint action, simultaneous on (the chains of)
$S^{2}$ (via the antipodal map) and on $\bb{Z}^{2}$ via $\rho$. It is easy to see that
$H^{0}(S^{2};\bb{Z}^{2}_{\rho})\approx \bb{Z}^{2}_{\rho}$ as $\bb{Z}[\bb{Z}_{2}]$-modules, and
$H^{2}(S^{2};\bb{Z}^{2}_{\rho})\approx \bb{Z}^{2}$, i.e., $\bb{Z}^{2}$ with the trivial
$\bb{Z}_{2}$-action.

A direct computation (e.g., see~\cite[pages~58-9]{kB}) yields the following values for
$E^{p,q}_{2}$:
\[
    E^{p,q}_{2}=
    \begin{cases}
      \bb{Z}^{2}    &\text{if $(p,q)=(0,2)$;}\\
      (\bb{Z}_{2})^{2}  &\text{if $q=0$ and $p$ odd, or if $q=2$ and $p>0$ and even;}\\
      0,        &\text{otherwise.}
    \end{cases}
\]
It follows easily from this that we have an exact sequence
\[0\rightarrow H^{2}(\bb{R}P^{2};\bb{Z}^{2}_{\rho})\overset{\pi^{\ast}}{\rightarrow}H^{2}(S^{2};\bb{Z}^{2}_{\rho})=\bb{Z}^{2}\rightarrow (\bb{Z}_{2})^{2}\rightarrow 0.\]
Thus, $H^{2}(\bb{R}P^{2};\bb{Z}^{2}_{\rho})\approx \bb{Z}^{2}$, and $\pi^{\ast}$ is injective.
Therefore, in this case Theorem~\ref{T:main} reduces to the following, which is an elaboration of
the second statement of Proposition~\ref{P:maingenuszero}(b):
\begin{prop}
The symplectic bundles $\xi: T^{2}\rightarrow E\rightarrow \bb{R}P^{2}$ inducing a non-trivial
$\bb{Z}_{2}$-module structure $\bb{Z}^{2}_{\rho}$ on $H_{1}(T^{2})=\bb{Z}^{2}$ are classified by
$H^{2}(\bb{R}P^{2};\bb{Z}^{2}_{\rho})\approx \bb{Z}^{2}$.  For such a $\xi$, E admits a closed
$2$-form satisfying (1) if and only if $c(\xi)=0$.
\end{prop}
\begin{proof}
The foregoing calculation implies the first statement of the proposition.  The second follows
from the injectivity of $\pi^{\ast}$, Lemma~\ref{L:pullback}, and
Proposition~\ref{P:maingenuszero}(a).
\end{proof}

\vspace{.2in}

This concludes our proof of Theorem~\ref{T:main}.

\vspace{.2in}

\section{Proofs of the main corollaries}~\label{S:Cor}

\vspace{.2in}

\noindent\emph{Proof of Corollary~\ref{C:finite}}:

For any connected surface $B$, $H^{2}(B;\bb{Z}^{2}_{\rho})$ is a finitely-generated abelian
group, hence, its torsion subgroup is finite.  The result now follows from
Proposition~\ref{P:classification} and Theorem~\ref{T:main}.  $\quad \square$

\vspace{.2in}

\noindent\emph{Proof of Corollary~\ref{C:principal}}:

The group of a principal torus bundle is $T^{2}$ acting on itself by translations.  If $\sigma$
denotes the standard symplectic form on $T^{2}$, then the translations clearly preserve $\sigma$,
i.e., $T^{2}\subseteq\ Symp(T^{2},\sigma)$, so the bundle has a canonical symplectic structure.
The corresponding representation
\[ \rho:\pi_{1}(B)\rightarrow\ \pi_{0}(Symp(T^{2},\sigma))\]
factors through $\pi_{0}(T^{2})=0$, so it is trivial.  Hence, when $B$ is a connected surface,
the only cases in which the characteristic classes $c(\xi)\in H^{2}(B;\mathbb{Z}^2)$ do not have
finite order are when $B$ is closed and orientable and $\xi$ is non-trivial. $\quad  \square$

\vspace{.2in}

\noindent\emph{Proof of Corollary~\ref{C:torusaction}}:

Let $i:T^{2}\rightarrow\ E$ be the inclusion onto a fixed orbit, and let $p:E\rightarrow\
E/T^{2}$ be the usual projection.  Then, it is a standard fact that $T^{2}\rightarrow\
E\rightarrow\ E/T^{2}$ is a principal torus bundle, say $\xi$. By Corollary~\ref{C:principal},
$\xi$ has a canonical structure as a symplectic torus bundle.  Let $\sigma$ be the standard
symplectic form on $T^{2}$, and let $\sigma_{b}$ be the corresponding symplectic forms on the
fibres (equiv., orbits).  By hypothesis, $E$ admits a symplectic form with respect to which all
the fibres are symplectic submanifolds. Thus, the restriction map
$i^{\ast}:H^{2}_{DR}(E)\rightarrow\ H^{2}_{DR}(T^{2})$ is surjective, and, by Thurston's result,
there is a closed $2$-form $\beta$ on $E$ satisfying condition (1), that is,
$\beta|T^{2}_{b}=\sigma_{b}$, for all $b\in E/T^{2}$. Assuming that the closed, connected surface
$E/T^{2}$ is orientable, we can then apply the preceding corollary to conclude that $\xi$ is
trivial, \emph{as a symplectic torus bundle}. Thus, it admits a section. But the existence of a
section is independent of the group of the bundle. Therefore, $\xi$ has a section as a principal
$T^{2}$ bundle, and, and therefore it is trivial as a principal $T^{2}$ bundle, which implies the
stated result.  It remains to verify that $E/T^{2}$ is orientable.  But this follows from  a
standard fact about smooth fibre bundles that are orientable, that is, for which the fibres can
be given orientations that are locally coherent over the base.  For such a bundle---for example
$\xi$--- the orientability of the base is equivalent to the orientability of the total space.
$\quad \square$

\appendix

\vspace{.3in}

\begin{center}
\textbf{APPENDIX}
\end{center}

The main arguments of the paper make use of certain known classification results in order to pass
from statements about smooth fibre bundles to statements about group extensions.  The following
three appendices briefly  explain these results, starting with facts about torus bundles, then
passing to the classification of $K(A,1)$-fibrations, and ending with a comparison between that
classification and the classification of corresponding group extensions.

\section{$T^{2}$-bundles and $T^{2}$-fibrations}\label{A:bundles}

Let $\mathcal{E}(T^{2})$  (resp., $\mathcal{E}_{+}(T^{2})$) denote the monoid of self homotopy
equivalences (resp., orientation-preserving self homotopy equivalences)  of $T^{2}$. These
receive the compact-open topology. Let  $Diff_{+}(T^{2})$ (resp.,$ Diff_{0}(T^{2})$) denote the
subgroup of orientation-preserving diffeomorphisms of $T^{2}$  (resp., the identity component of
$Diff(T^{2}))$.  Finally, let $\omega$ be any symplectic form on $T^{2}$, and let
$Symp(T^{2},\omega)$ be the group of symplectomorphisms of $(T^{2},\omega)$. These groups of
diffeomorphisms are usually given the $C^{k}$ topology, for $1\leq k\leq\infty$.  The choice of
$k$ does not make a difference for our discussion.  Regarding $T^{2}$ as acting on itself by
translation, we have $T^{2}\subseteq Diff_{0}(T^{2})$.

\begin{prop}\label{P:torusdiff}
The following inclusions are homotopy equivalences:
\begin{enumerate}
\item $T^{2}\rightarrow Diff_{0}(T^{2})$.
\item $Diff(T^{2})\rightarrow \mathcal{E}(T^{2})$.
\item $ Diff_{+}(T^{2})\rightarrow \mathcal{E}_{+}(T^{2})$.
\item $Symp(T^{2},\omega)\rightarrow  Diff_{+}(T^{2})$
\end{enumerate}
\end{prop}
\begin{proof}

(a),(b): These are well-known results, due originally to Earle and Eells (cf.,
Gramain~\cite{aG}). (c) follows immediately from (b). (d): Given any orientation-preserving
diffeomorphism $h$, $h^{\ast}(\omega)$ is homologous to $\omega$, since $h$ has degree one. Thus,
Moser's method (\cite[pages 93-97]{MS}) may be applied to the family of symplectic forms
$\omega_{t}=(1-t)\omega+th^{\ast}\omega$,  producing an isotopy $\psi_{t}$ between the identity
and a diffeomorphism $\psi_{1}$ that satisfies $\psi_{1}^{\ast}h^{\ast}\omega=\omega$. Therefore,
$h_{t}=h\psi_{t}$ is an isotopy between $h$ and a symplectomorphism $h_{1}$.  The isotopy can be
constructed so as to be continuous in $h$ and  remain in $Symp$ if $h$ is a symplectomorphism. It
follows that the map given by $h\mapsto h_{1}$ is a homotopy inverse for the inclusion map.
\end{proof}

Since, as is well known, $B\mathcal{E}_{+}(T^{2})$ classifies oriented $T^{2}$-fibrations,
statements (c) and (d) immediately give the following result.

\begin{corr}
Let $B$ be a smooth, connected manifold.  Equivalence classes of symplectic torus bundles over
$B$ correspond bijectively to fibre-homotopy equivalence classes of oriented $T^{2}$-fibrations
over $B$.$\quad \square$
\end{corr}

\section{K(A,1)-fibrations}\label{A:fibrations}

Let $A$ be an abelian group.  Following C.A. Robinson~\cite{cR}, for any $n\geq 1$, we let
$K(A,n)$ be a CW complex which is a topological abelian group on which $Aut(A)$ acts by cellular
automorphisms.  Let $Q$ be a CW complex of type $K(Aut(A),1)$ and $\tilde{Q}$ its universal
cover. Thus, there is a free, diagonal left-action of $Aut(A)$ on the cartesian product
$K(A,2)\times \tilde{Q}$, with respect to which the projection
$K(A,2)\times\tilde{Q}\rightarrow\tilde{Q}$ is equivariant.  Therefore, it descends to a
fibration $p:\hat{K}(A,2)\rightarrow Q$ with fibre $K(A,2)$.

Robinson shows that  $\hat{K}(A,2)$ classifies Hurewicz fibrations with fibres of the homotopy
type of $K(A,1)$ and base spaces of the homotopy type of a CW complex.  Thus, over such a base
space $B$, the fibre-homotopy equivalence classes of $K(A,1)$- fibrations correspond bijectively
to homotopy classes of maps $B\rightarrow \hat{K}(A,2)$. Throughout this paper, we use the
`based' convention for equivalences, whereby each base space has a basepoint and each fibre has a
fixed identification with a given space. See~\cite{cR} and~\cite[16.7]{aD}.

\begin{remark}
 By Proposition~\ref{P:torusdiff},  $BDiff(T^{2})$  classifies $T^{2}$-fibrations, which are the same as $K(\bb{Z}^{2},1)$-fibrations.  So $BDiff(T^{2})$ is homotopy equivalent to $\hat{K}(\bb{Z}^{2},2)$, implying that it too fibres over
over $K(GL(2,\bb{Z}^{2}),1)$ with fibre $K(\bb{Z}^{2},2)$. This fact is well known, but we
mention it to connect the two constructions of classifying spaces. It gives one way of seeing
why, for a simply-connected base space, there is a bijective correspondence between equivalence
classes of torus bundles and equivalence classes of principal torus bundles. A similar comment
applies to $BSymp(T^{2},\omega)$, which fibres over $K(SL(2,\bb{Z}^{2}),1)$ with fibre
$K(\bb{Z}^{2},2)$
\end{remark}

 As usual, each $K(A,1)$-fibration admits a representation $\rho:\pi_{1}(B)\rightarrow Aut(A)=\pi_{0}(\mathcal{E}(K(A,1)))$. We are interested in the finer classification that fixes such a $\rho$.  Robinson derives this from his construction as follows.
The fibration $p:\hat{K}(A,2)\rightarrow Q$ admits a canonical section $s_{0}:Q\rightarrow
\hat{K}(A,2)$ defined by the rule $s_{0}[q]=[\vartheta,q]$, where here $[\quad]$ refers to the
$Aut(A)$-orbit and $\vartheta$ denotes the identity element of the abelian group $K(A,2)$
Clearly, representations $\rho:\pi_{1}(B)\rightarrow Aut(A)$ correspond to homotopy classes of
maps $r:B\rightarrow Q$, whereas $K(A,1)$-fibrations over $B$ with associated representation
$\rho$ correspond to homotopy classes of maps $f:B\rightarrow \hat{K}(A,2)$ for which $pf$
induces $\rho$. In fact, as Robinson shows, if we fix $\rho$ ( and $r$ inducing $\rho$),  the
foregoing set of homotopy classes may be described as the set of homotopy classes of lifts $f$ of
$r$ to $\hat{K}(A,2)$. Let $f_{0}$ denote the lift $s_{0}r$.

Given two lifts $f$ and $g$ of $r$, classical obstruction theory produces a so-called primary
obstruction class $d(f,g)\in H^{2}(B;\pi_{2}(K(A,2)))=H^{2}(B;A_{\rho})$ whose vanishing is a
necessary condition for the existence of a homotopy of lifts between $f$ and $g$. In this
context, the condition is also sufficient. Moreover, given any $g$ and any class $d\in
H^{2}(B;A_{\rho})$, there is a unique homotopy class of lifts $f$ such that $d(f,g)=d$.   We now
set $d(f,f_{0})=c(f)$, where $f_{0}$ is the lift described above.  A standard additivity formula
yields $d(f,g)=c(f)-c(g)$.  If $f$ classifies a $K(A,1)$-fibration $\eta$, we may write
$c(f)=c(\eta)$. This is the so-called characteristic class of $\eta$ that we have been using. It
follows that the rule $\eta\mapsto c(\eta)$ gives a bijection between equivalence classes of
fibrations and $H^{2}(B, A_{\rho})$, as stated earlier.

There remains one further observation about the classes $c(\eta)$ which we have used in an
important way. Start by considering a homotopy $h_{t}$ between the lifts $f$ and $f_{0}$ of $r$
as above. Then, for any $b\in B$, $h_{t}(b)$ is a path from $f(b)$ to $f_{0}(b)$ lying completely
in a fibre of $p:\hat{K}(A,2)\rightarrow Q$.  This motivates the following construction given by
Robinson: Let $P$ denote the space of all paths $\gamma$ in $\hat{K}(A,2)$ that begin in
$s_{0}(Q)$ and lie completely in a fibre of $p$. The rule $\gamma\mapsto \gamma(1)$ defines a
fibration $P\rightarrow \hat{K}(A,2)$ with fibre of type $K(A,1)$.  Clearly the partial
homotopies between any lift $f$ of $r$ and the lift $f_{0}$ correspond to partial lifts of $f$ to
$P$.  It follows easily that we can interpret $c(f)$ as the primary obstruction to lifting $f$ to
$P$.

Now, Robinson shows that $P$ is, in fact, a \emph{universal} $K(A,1)$-fibration over
$\hat{K}(A,2)$, so that $f^{\ast}(P)$ is equivalent to $\eta$.  This implies that $c(\eta)$ may
be interpreted directly as the primary obstruction to a section of $\eta$, which is the
interpretation we have used.

\section{Extensions by A}

Let $\bb{S}: A\rightarrowtail G\twoheadrightarrow \pi$ be an extension of a group $\pi$ by the
abelian group $A$.  There is a corresponding $K(A,1)$-fibration, which we write as $\eta:
K(A,1)\overset{i}{\rightarrow} K(G,1)\overset{p}{\rightarrow} K(\pi,1)$.  Of course, the homotopy
exact sequence of $\eta$ collapses to $\bb{S}$.  We use $i_{\ast}$ and $p_{\ast}$ to denote the
corresponding homomorphisms in $\bb{S}$. The representation $\rho$ corresponding to $\eta$ is the
same as that induced by inner automorphisms of $G$ in $\bb{S}$. Let us hold this fixed.

Let $f:\pi\times\pi\rightarrow A$ be the normalized $2$-cocycle of $\bb{S}$. In terms of the bar
resolution of $\pi$, we may write $f$ as the (possibly infinite) formal sum $\Sigma f(x,y)[x|y]$,
where $x,y$ range over $\pi$.  In this appendix we show how this sum can be recognized as the
primary obstruction to sectioning $\eta$.  This establishes the identification
$c(\bb{S})=c(\eta)$, which we have been using throughout the paper.  This fact is certainly part
of the classical folklore of the subject, but I have been unable to find an explicit reference.

The description of the primary obstruction  can be conveniently simplified in this case by using
the following observation, which follows almost immediately from definitions.

\begin{lemma}
Let $\zeta: F \overset{i}{\rightarrow} E \overset{p}{\rightarrow} B$ be a fibration, with $F$
connected and $B$ a connected CW complex, and assume that
$i_{\ast}:\pi_{m-1}(F)\rightarrow\pi_{m-1}(E)$ is injective.  Let $\sigma:B^{m-1}\rightarrow E$
be a section of $\zeta$ over the $m-1$-skeleton of $B$, and let $o(\sigma)$ denote the
obstruction cocycle to extending $\sigma$ over the $m$-skeleton.    Finally, suppose that if
$\chi:D^{m}\rightarrow B$ is the characteristic map of an $m$-cell $e$ of $B$, then
$\sigma\chi|\partial D^{m}:\partial D^{m}\rightarrow E$ represents $i_{\ast}(\alpha)\in
\pi_{m-1}(E)$.  Then
\[o(\sigma)(e)=\alpha.\quad\quad\square \]
\end{lemma}

The best framework for recognizing $\Sigma f(x,y)[x|y]$ as the desired obstruction cocycle is
that of semisimplicial topology, as in~(\cite[Chapters 1--3]{kL}).  Thus, for example, we can
describe $K(\pi,1)$ semisimplicially as consisting of one $0$-simplex, denoted $[\quad]$, and a
$k$-simplex for each integer $k\geq 1$ and each symbol $[x_{1}|\ldots|x_{k}]$, where
$x_{1},\ldots,x_{k}$ range over $\pi$, with the well-known face and degeneracy maps. Similarly
for $K(G,1)$.  The surjection $p_{\ast}:G\twoheadrightarrow \pi$ shows how to map $K(G,1)$ onto
$K(\pi,1)$.  This map is a minimal Kan fibration, say $\kappa$~\cite[page 64]{kL}. We shall
define an obstruction to sectioning $\kappa$.

Let $s:\pi\rightarrow G$ be a function that is a right inverse of $p_{\ast}$ and is related to
the $2$-cocycle $f:\pi\times\pi\rightarrow\ G$ by equation~(\ref{E:cocycle}).  Use $s$ to define
a (semisimplicial) section $\sigma$ of $\kappa$ over the $1$-skeleton of $K(\pi,1)$: this is
determined by $\sigma[\quad]=[\quad]$ and $\sigma[x]=[s(x)]$.  Note that each $1$-simplex
$[s(x)]$ determines a directed loop in $K(G,1)$, say $<s(x)>$; these may be concatenated. Now
consider any $2$-simplex $[x|y]$ of $K(\pi,1)$. Its boundary consists of the $1$-simplexes
$\partial_{0}[x|y]=[y]$, $\partial_{1}[x|y]=[xy]$, and $\partial_{2}[x|y]=[x]$, with
corresponding loops concatenated as $<x><y><xy>^{-1}$.  Therefore, the loop obtained by applying
$\sigma$ to the boundary is $<s(x)><s(y)><s(xy)>^{-1}$.  Using the semisimplicial homotopy
relation in $K(G,1)$ and the definition of $f(x,y)$, this loop is easily shown to be homotopic to
$<i_{\ast}(f(x,y))>$ in $K(G,1)$. Thus, it represents $i_{\ast}(f(x,y))$. It now follows from the
semisimplicial analog of the above lemma that the obstruction to extending $\sigma$ over the
$2$-skeleton is precisely the cocycle $\Sigma f(x,y)[x|y]$, as desired.

The foregoing can be translated to the more conventional topological obstruction theory by
applying the geometric realization functor. This transforms $\kappa$ into a topological fibration
equivalent to $\eta$ and $\sigma$ into a partial section producing the same obstruction. Thus
$c(\eta)=c(\bb{S})$.

\end{document}